\documentclass[12pt,reqno]{amsart}
\usepackage{enumerate, latexsym, amsmath, amsfonts, amssymb, amsthm, color}
\def\pmod #1{\ ({\rm{mod}}\ #1)}
\def\Z{\Bbb Z}
\def\N{\Bbb N}

\def\l{\left}
\def\r{\right}
\def\bg{\bigg}
\def\({\bg(}
\def\){\bg)}
\def\t{\text}
\def\f{\frac}
\def\mo{{\rm{mod}\ }}

\def\ls{\leqslant}
\def\gs{\geqslant}

\def\bi{\binom}

\def\eq{\equiv}

\def\Proof{\noindent{\it Proof}}

\theoremstyle{plain}
\newtheorem{theorem}{Theorem}

\newtheorem{lemma}{Lemma}

\theoremstyle{definition}

\theoremstyle{remark}
\newtheorem{remark}{Remark}

 \vspace{4mm}

\begin{document}

\hbox{Preprint, {\tt arXiv:2006.16089}}
\medskip

\title
[{A new extension of the Sun-Zagier result}]
{A new extension of the Sun-Zagier result involving Bell numbers and derangement numbers}

\author
[Zhi-Wei Sun] {Zhi-Wei Sun}

\address{Department of Mathematics, Nanjing
University, Nanjing 210093, People's Republic of China}
\email{zwsun@nju.edu.cn}

\keywords{Bell number, Bell polynomial, congruence, derangement number.
\newline \indent Supported by the Natural Science Foundation of China (grant no. 11971222).}

\begin{abstract}  Let $p$ be any prime and let $a$ and $n$ be positive integers with $p\nmid n$. We show that $$\sum_{k=1}^{p^a-1}\frac{B_k}{(-n)^k}\equiv a(-1)^{n-1}D_{n-1}\pmod {p},$$
where $B_0,B_1,\ldots$ are the Bell numbers and $D_0,D_1,\ldots$ are the derangement numbers.
This extends a result of Sun and Zagier published in 2011. Furthermore, we prove that
$$(-x)^n\sum_{k=1}^{p^a-1}\frac{B_k(x)}{(-n)^k}\equiv -\sum_{r=1}^ax^{p^r}\sum_{k=0}^{n-1}\frac{(n-1)!}{k!}(-x)^k\pmod{p\mathbb Z_p[x]},$$
where $B_k(x)=\sum_{l=0}^kS(k,l)x^l$ is the Bell polynomial of degree $k$ with $S(k,l)\ (0\ls l\ls k)$ the Stirling numbers of the second kind,
and $\mathbb Z_p$ is the ring of all $p$-adic integers.
\end{abstract}
\maketitle

\section{Introduction}
\setcounter{lemma}{0}
\setcounter{theorem}{0}
\setcounter{corollary}{0}
\setcounter{remark}{0}
\setcounter{equation}{0}

Let $B_0=1$. For each $n\in\Z^+=\{1,2,3,\ldots\}$ let $B_n$ denote the number of
  partitions of a set of cardinality $n$.  For example, $B_3=5$ since there are totally 5 partitions of $\{1,2,3\}$:
  $$\{\{1\},\{2\},\{3\}\},\{\{1,2\},\{3\}\},\ \{\{1,3\},\{2\}\},\ \{\{2,3\},\{1\}\},\{\{1,2,3\}\}.$$
The Bell numbers $B_0,B_1,\ldots$, named after E. T. Bell who studied them in the 1930s,
play important roles in combinatorics.
Here are values of $B_1,\ldots,B_{7}$:
$$B_1=1, \ B_2=2,\ B_3=5,\ B_4=15,\ B_5=52,\ B_6=203,\ B_7=877.$$
It is known that
$$\sum_{n=0}^\infty B_n\f{x^n}{n!}=e^{e^x-1}\  \t{and}\ \ B_{n+1}=\sum_{k=0}^n\bi nk B_k\ (n=0,1,2,\ldots).$$
The author's conjecture (cf. \cite[Conjecture 3.2]{S11}) that the sequence $(\root{n+1}\of{B_{n+1}}/\root n\of{B_n})_{n\gs1}$
is strictly decreasing (with limit $1$), is still open.
For any prime $p$ and  $m,n\in\N=\{0,1,2,\ldots\}$, we have the classical Touchard congruence (cf. \cite{T})
$$B_{p^m+n}\eq mB_n+B_{n+1}\pmod p.$$

Let $D_0=1$, and define $D_n$ $(n\in\Z^+)$ by
$$D_n=|\{\pi\in S_n:\ \pi(k)\not=k\ \t{for all}\ k=1,\ldots,n\}|.$$
Those $D_0,D_1,D_2,\ldots$ are called the
derangement numbers, and they were first introduced by Euler. It is well known that
$$D_n=n!\sum_{k=0}^n\f{(-1)^k}{k!}\ \ \t{for all}\ n\in\N.$$

In 2011, the author and D. Zagier \cite{SZ} showed that for any prime $p$ and $n\in\Z^+$ with $p\nmid n$ we have
\begin{equation}\label{SZ}\sum_{k=1}^{p-1}\f{B_k}{(-n)^k}\equiv (-1)^{n-1}D_{n-1}\pmod {p},
\end{equation}
which relates the Bell numbers to the derangement numbers.
The surprising congruence \eqref{SZ} was called the Sun-Zagier congruence
by Y. Sun, X. Wu and J. Zhuang \cite{SWZ} who used the umbral calculus to give a generalization,
by I. Mez\H o and T. L. Ramirez \cite{MR} in 2017 who extended it to the so0called $r$-Bell numbers,
and by Q. Mu \cite{Mu} in 2018 who re-proved via an identity of R.J. Clarke and M. Sved \cite{CS} relating the Bell numbers to the derangement numbers.

In this paper we extend the fundamental Sun-Zagier result in a new way.

\begin{theorem}\label{Th1.1}
 Let $p$ be any prime and let $a$ be a positive integer. For any $n\in\Z^+$ with $p\nmid n$, we have
\begin{equation}\label{BD}\sum_{k=1}^{p^a-1}\f{B_k}{(-n)^k}\equiv a(-1)^{n-1}D_{n-1}\pmod {p}.
\end{equation}
\end{theorem}
\begin{remark}\label{Rem1.1} Note that \eqref{BD} in the case $a=1$ gives \eqref{SZ}.
\end{remark}

For $n\in\Z^+$ and $k\in\{0,\ldots,n\}$, the Stirling number $S(n,k)$ of the second kind
denotes the number of ways to partition the set $\{1,\ldots,n\}$ into $k$ disjoint nonempty parts.
In addition, we adopt the usual convention $S(0,0)=1$. For $n\gs k\gs0$, it is well known that
\begin{equation}\label{Stirling}k!S(n,k)=\sum_{j=0}^k\bi kj(-1)^{k-j}j^n.
\end{equation}

For any $n\in\N$, the Bell polynomial (or the Touchard polynomial) of degree $n$ is given by
\begin{equation}\label{BellP}B_n(x)=\sum_{k=0}^n S(n,k)x^k.
\end{equation}
Clearly, $B_n(1)=B_n$ for all $n\in\N$, and $B_n(x)=x\sum_{k=1}^nS(n,k)x^{k-1}$ for all $n\in\Z^+$.
Theorem \ref{Th1.1} actually follows from our following theorem concerning the Bell polynomials.

\begin{theorem}\label{Th1.2} Let $a$ be any positive integer.
For any $n\in\Z^+$ and prime $p\nmid n$, we have
\begin{equation}\label{BDx}(-x)^n\sum_{k=1}^{p^a-1}\f{B_k(x)}{(-n)^k}\eq -\sum_{r=1}^ax^{p^r}\sum_{k=0}^{n-1}\f{(n-1)!}{k!}(-x)^k\pmod{p\Z_p[x]},
\end{equation}
where $\Z_p$ denotes the ring of all $p$-adic integers.
\end{theorem}
\begin{remark}\label{Th1.2} The congruence \eqref{BDx} in the case $a=1$
was deduced by Sun and Zagier \cite{SZ} via the usual explicit formula \eqref{Stirling} for Stirling numbers of the second kind. Our Theorem 1.2 can be further extended in the spirit of \cite{SWZ,MR}, we omit the details.
\end{remark}

We will show Theorem 1.2 in the next section.

\section{Proof of Theorem 1.2}
\setcounter{lemma}{0}
\setcounter{theorem}{0}
\setcounter{corollary}{0}
\setcounter{remark}{0}
\setcounter{equation}{0}

\begin{lemma}\label{Lem2.1} Let $p$ be a prime and let $a\in\Z^+$.

{\rm (i)} For any $j,k\in\N$ with $j+k\ls p^a-1$, we have
$$\bi{p^a-1-k}j\bg/\bi{-1-k}j\eq1\pmod p.$$
In particular,
$$\bi{p^a-1}j\eq(-1)^j\pmod p\ \ \ \t{for all}\ j=0,\ldots,p^a-1.$$

{\rm (ii)} We have
$$B_{p^a}(x)\eq \sum_{r=0}^a x^{p^r}\pmod {p\Z_p[x]}.$$
\end{lemma}
\Proof. (i) Since $j+k<p^a$ we have
$$\f{\bi{p^a-1-k}j}{\bi{-1-k}j}=\f{\prod_{0<i\ls j}\f{p^a-i-k}i}{\prod_{0<i\ls j}\f{-i-k}i}
=\prod_{0<i\ls j}\l(1-\f{p^a}{i+k}\r)\eq1\pmod p.$$
When $k=0$, this yields
$$\bi{p^a-1}j\eq\bi{-1}j=(-1)^j\pmod p.$$

(ii) By A. Gertsch and A. M. Robert's extension \cite{GR} of Touchard's  congruence, for any $n\in\N$ we have
$$B_{p^a+n}(x)\eq B_{n+1}(x)+B_n(x)\sum_{r=1}^ax^{p^r}\pmod {p\Z_p[x]}.$$
In particular,
$$B_{p^a}(x)\eq B_1(x)+B_0(x)\sum_{r=1}^ax^{p^r}=x+\sum_{r=1}^ax^{p^r}=\sum_{r=0}^a x^{p^r}\pmod {p\Z_p[x]}.$$

In view of the above, we have completed the proof of Lemma \ref{Lem2.1}. \qed

\medskip
\noindent{\it Proof of Theorem \ref{Th1.2}}. It is known that
\begin{equation}\label{recur}B_{m+1}(x)=x\sum_{k=0}^m\bi mkB_k(x)\ \ \t{for all}\ m\in\N.
\end{equation}
In light of this and Lemma 2.1, for any prime $p$ we have
\begin{align*}\sum_{k=1}^{p^a-1}(-1)^kB_k(x)\eq&\sum_{k=1}^{p^a-1}\bi{p^a-1}kB_k(x)=\f{B_{p^a}(x)}x-B_0(x)
\\\eq& \sum_{r=1}^a x^{p^r-1}\pmod{p\Z_p[x]}.
\end{align*}
So the desired result holds when $n=1$.

Now we fix $n\in\Z^+$ and assume that \eqref{BDx} holds
 for every prime $p\nmid n$.

  Let $p$ be any prime not dividing $n+1$.
If $p\mid n$, then $n!/k!\eq0\pmod p$ for all $k=0,\ldots,n-1$, and hence
\begin{align*}(-x)^{n+1}\sum_{k=1}^{p^a-1}\f{B_k(x)}{(-n-1)^k}\eq&(-x)^{n+1}\sum_{k=1}^{p^a-1}\f{B_k(x)}{(-1)^k}
\eq (-x)^{n+1}\sum_{r=1}^ax^{p^r-1}
\\\eq& -\sum_{r=1}^ax^{p^r}\sum_{k=0}^n\f{n!}{k!}(-x)^k\pmod{p\Z_p[x]}.
\end{align*}

Now we suppose that $p\nmid n$. In view of \eqref{recur} and Lemma \ref{Lem2.1}, we have
\begin{align*} \sum_{k=1}^{p^a-1}\f{B_k(x)/x}{(-n)^k}=&\sum_{k=1}^{p^a-1}\f{\sum_{l=0}^{k-1}\bi{k-1}lB_l(x)}{(-n)^k}
=\sum_{l=0}^{p^a-2}\f{B_l(x)}{(-n)^l}\sum_{k=l+1}^{p^a-1}\f{\bi{k-1}l}{(-n)^{k-l}}
\\=&\sum_{l=0}^{p^a-2}\f{B_l(x)}{(-n)^{l+1}}\sum_{r=1}^{p^a-1-l}\f{\bi{l+r-1}{r-1}}{(-n)^{r-1}}
\\=&\sum_{l=0}^{p^a-2}\f{B_l(x)}{(-n)^{l+1}}\sum_{r=1}^{p^a-1-l}\f{\bi{-l-1}{r-1}}{n^{r-1}}
\\\eq&\sum_{l=0}^{p^a-2}\f{B_l(x)}{(-n)^{l+1}}\sum_{r=1}^{p^a-1-l}\bi{p^a-1-l}{r-1}n^{-(r-1)}
\\\eq&\sum_{l=0}^{p^a-1}\f{B_l(x)}{(-n)^{l+1}}\(\l(1+\f1n\r)^{p^a-1-l}-\f1{n^{p^a-1-l}}\)
\\\eq&\sum_{l=0}^{p^a-1}\f{n^lB_l(x)}{(-n)^{l+1}}\l(\f1{(n+1)^l}-1\r)\ (\mo\ p)
\end{align*}
with the aid of Fermat's little theorem.
Therefore
$$-n\sum_{k=1}^{p^a-1}\f{B_{k}(x)/x}{(-n)^{k}}\eq\sum_{k=1}^{p^a-1}\f{B_k(x)}{(-n-1)^{k}}
-\sum_{l=1}^{p^a-1}\f{B_l(x)}{(-1)^l}
\pmod {p\Z_p[x]}$$
and hence
$$\sum_{k=1}^{p^a-1}\f{B_k(x)}{(-n-1)^k}\eq-n\sum_{k=1}^{p^a-1}\f{B_k(x)/x}{(-n)^k}
+ \sum_{r=1}^ax^{p^r-1}\pmod{p\Z_p[x]}.$$
Combining this with \eqref{BDx}, we obtain
\begin{align*}(-x)^{n+1}\sum_{k=1}^{p^a-1}\f{B_k(x)}{(-n-1)^k}
\eq&-n\sum_{r=1}^ax^{p^r}\sum_{k=0}^{n-1}\f{(n-1)!}{k!}(-x)^k
+(-x)^{n+1}\sum_{r=1}^ax^{p^r-1}
\\=&-\sum_{r=1}^ax^{p^r}\sum_{k=0}^n\f{n!}{k!}(-x)^k\pmod{p\Z_p[x]}.
\end{align*}
This concludes the induction step.

By the above, the proof of Theorem 1.1 is now complete. \qed

\end{document}